%% file: main.tex
\documentclass{article}
\usepackage[utf8]{inputenc}
\usepackage[centertags]{amsmath}
\usepackage{amsfonts,amssymb,amsthm}
\usepackage{float}
\usepackage{graphicx,subcaption,url}
\usepackage{tikz}
\usepackage{bussproofs}

\usepackage{lineno}

\newenvironment{scprooftree}[1]%
  {\gdef\scalefactor{#1}\begin{center}\proofSkipAmount \leavevmode}%
  {\scalebox{\scalefactor}{\DisplayProof}\proofSkipAmount \end{center} }

\title{Reflecting on beauty: the aesthetics of mathematical discovery}

\author{
Filip D. Jevti\'{c} \\{\small Mathematical Institute}\\[-2mm] {\small SASA,  Belgrade}
\and
Jovana Kosti\'{c} \\{\small Faculty of Philosophy}\\[-2mm] {\small University of Belgrade}
\and
Katarina Maksimovi\'{c} \\{\small Faculty of Philosophy}\\[-2mm] {\small University of Belgrade}
}
\date{}

\begin{document}

\maketitle

\begin{abstract}
Mathematical research is often motivated by the desire to reach a beautiful result or to prove it in an elegant way. Mathematician's work is thus strongly influenced by his aesthetic judgments. However, the criteria these judgments are based on remain unclear. In this article, we focus on the concept of mathematical beauty, as one of the central aesthetic concepts in mathematics. We argue that beauty in mathematics reveals connections between apparently non-related problems or areas and allows a better and wider insight into mathematical reality as a whole. We also explain the close relationship between beauty and other important notions such as depth, elegance, simplicity, fruitfulness, and others. 
 
\end{abstract}


\begin{quote}
\textit{In memory of Kosta Do\v{s}en on the occasion of his seventieth birthday.}
\end{quote}

\section{Introduction}
There are very few things in the world as widely used and as indispensable in both science and everyday life as mathematics. However, the motivation for mathematical research generally does not lie in its utility. Mathematicians do not tend to choose the subject of their study because of its projected applications. Rather, they get intrigued by the subject itself, their research often being motivated by the desire to reach a \textit{beautiful} result or to prove it in an \textit{elegant} way. As Rota observed, ``The lack of beauty in a piece of mathematics is of frequent occurrence, and it is a strong motivation for further mathematical research'' (\cite{rota}, p. 178). Mathematician's work is thus strongly influenced by his aesthetic judgments.\footnote{About the significance of aesthetic judgments for science see \cite{mcallister1}.} On the other hand, the main objective of mathematical research is truth. Above all, we want our theories to be consistent and our results to be correct. It is not immediately clear how these two objectives are connected and why aesthetic properties such as beauty serve as the driving force for something that is ultimately unrelated to aesthetics. 

In this article, we analyze the concept of mathematical beauty, as one of the central aesthetic concepts in mathematics and examine its relationship with non-aesthetic properties of mathematical results, such as truth, simplicity, fruitfulness, etc., as well as other aesthetic concepts such as elegance and depth. We don't strive to offer a complete characterization of all aesthetic aspects of mathematics (for a discussion on some of those aspects not mentioned here see \cite{thomas}), nor do we aim to give a psychological or a phenomenological account of aesthetic experience. But we do hope that our discussion will shed some light on the nature of mathematics itself and what we believe to be some of the key problems to understanding the subject of beauty and the aesthetics of mathematical discovery.  

Some authors denied the possibility of making any proper aesthetic judgments involving beauty in mathematics. They assumed beauty must have a perceptual component which mathematics lacks \cite{gerwen}. It is easy to see that this argument is unfounded, as non-perceptual beauty is present not only in mathematics but also in other domains, such as literary art (\cite{cellucci}, section 7).  So, the first step in our analysis of mathematical beauty would be to acknowledge its conceptual nature and the fact that this specific type of beauty requires prior knowledge to be recognized. Mathematical beauty is ``dependent upon the mathematician’s
background knowledge, so it is usually recognizable only to the well trained'' (\cite{cellucci} p. 346). Yet, there is still a lot to be done if we are to fully understand this concept.  

In the literature addressing the question of beauty in mathematics, two main strands can be distinguished. 
The first one suggests that the aesthetic criteria in mathematics are subjective, reflecting the preferences and expertise of the individual making the aesthetic judgment, as well as the context in which the judgment occurs. According to this view, asserting the beauty of a mathematical object or a mathematical result does not derive its justification from an objective fact, but from the preferences of the subject making the assertion (e.g. \cite{rota}). 
The second strand takes the criteria for mathematical beauty as objective and grounded in the content judged as beautiful. According to this stance, anyone with sufficient knowledge of the content should be capable of recognizing its beauty. Theories aligned with this perspective often link mathematical beauty to epistemic qualities such as truth, clarity, explanatory power, depth, and fruitfulness. Some of them have even attempted to reduce mathematical beauty to these non-aesthetic qualities, claiming that aesthetic judgments in mathematics are not really aesthetic at all (\cite{todd}, also see \cite{harre}). 

The main problem with the first perspective is its inability to explain the connection between beauty and truth (and other epistemic qualities). ``Real beauty is for mathematicians inseparable from truth, but, more than that, it is a sign of truth, and moreover of important truth'' (\cite{DA}, p. 34). If beauty is relative to the subjective preferences of mathematicians, how can it serve not only as the intrinsic motivational force behind mathematical progress but also as the standard by which we judge the extent of this progress?\footnote{Here, we make a reasonable assumption that both mathematical progress and progress in science generally are based on objective facts rather than preferences.} On the other hand, the second view faces the problem of how to define the properties that constitute mathematical beauty while at the same time accounting for their specific \textit{aesthetic} qualities. Claiming that mathematical beauty is reducible to epistemic qualities might suggest that mathematicians, when making aesthetic judgments of their subject, misconceive its epistemic properties as aesthetic ones. In our view, this shows there is something deeply problematic with this claim as it seems to presuppose that mathematicians aren't competent enough to differentiate between their own aesthetic and epistemic judgments, while we possess better insight into the intentions of their statements. 
 
The account presented in this paper belongs to the second strand. We argue that mathematical beauty is objective, forming a distinct aesthetic category irreducible to non-aesthetic or other aesthetic categories from outside mathematics. However, we maintain that mathematical beauty is closely related to truth and other epistemic qualities, because it serves as a guiding light, directing our attention to what is fruitful and deep and guiding us towards important discoveries. Just like the beautiful color, shape, and overall aesthetically pleasing appearance of fruit helps us instantly spot pieces of fruit that are fresher, nutritionally richer, and free from harmful pathogens, beauty in mathematics helps us recognize those structures and properties of mathematical objects that are in the mathematical sense more fundamental. An advantage of our account lies in its ability to explain the enduring focus on symmetry, simplicity, and harmony within discussions of mathematical beauty, while also illustrating how these attributes are intertwined with the notion of truth in mathematics.

The structure of this paper is as follows. In the second section, we present the proposed account in more detail. As the perspective presented in this paper draws inspiration from Plato and his idea about the role of beauty in acquiring knowledge of Forms, in section 3 we briefly discuss Plato's views on beauty and their connection to mathematics. We then move on to discuss the three levels of mathematical entities: objects, theorems, and proofs, and the three corresponding notions of beauty in section 4. In section 5 we propose a unified account of mathematical beauty and discuss some possible problems and solutions. 

\section{A purpose-driven beauty}

The account we are proposing includes three components. As we already mentioned in the introduction, we claim that:
\begin{itemize}
    \item[$\cdot$] Beauty in mathematics is objective;
    \item[$\cdot$] It is not reducible to non-aesthetic epistemic properties;
    \item[$\cdot$] It is connected to truth and other epistemic properties by pointing out to and guiding us toward them.
\end{itemize}

Each of these claims requires some explanation. It might be best to begin with the following example given by the physicist Steven Weinberg (\cite{vajnberg}, p. 133):

``A physicist who says that a theory is beautiful does not mean quite the same thing that would be meant in saying that a particular painting or a piece of music or poetry is beautiful. It is not merely a personal expression of aesthetic pleasure; it is much closer to what a horse trainer means when he looks at a racehorse and says that it is a beautiful horse. The horse trainer is of course expressing a personal opinion, but it is an opinion about an objective fact: that, on the basis of judgments that the trainer could not easily put into words, this is the kind of horse that wins races. [...] The physicist’s sense of beauty is also supposed to serve a purpose – it is supposed to help the physicist select ideas that help us to explain nature.'' 

This example shows that we are indeed capable of contemplating and assessing beauty in at least two different ways. On one hand, the judgments about the beauty of nature are devoid of any purpose or context. We assess the beauty of an object solely based on how it affects us. For instance, we might call a running horse beautiful because witnessing a horse running free in the wild is simply a beautiful sight to behold and an aesthetically pleasing experience. 

On the other hand, there is a different kind of aesthetic judgment where beauty is assessed with a specific purpose in view. In this case, there is an activity or a predetermined goal, and objects are judged as beautiful to the extent to which they seem likely to fulfill that purpose or help us attain that goal. We do not admire the object itself but rather what it represents for that particular activity. In Weinberg’s example, such would be the trainer’s assertion that the horse is beautiful. The fact that we can clearly distinguish between these two contexts and acknowledge that a horse may be beautiful in the first sense but not in the second one (for example, having a long mane or a delicate neck may contribute to an aesthetically pleasing appearance, but a horse with these features might not be beautiful in the eyes of the horse trainer) shows that there are indeed two distinct types of aesthetic judgments. The distinction between the two types of aesthetics judgments is somewhat similar to Harr\'{e}'s distinction between first order and second order aesthetic appraisals (see \cite{harre}). 

But unlike Harr\'{e}, we believe both types of judgments are aesthetic because they possess a normative component and they are connected to a specific type of aesthetic experience, such as satisfaction or admiration for the properties that evoke the judgment. This enables one to distinguish an aesthetic judgment from a mere description of an object.\footnote{The inability to differentiate between normative and purely descriptive judgments is in philosophy called the naturalistic fallacy, which is a concept famously introduced by George Edward Moore in \cite{mur}.}

While Weinberg addresses the question of beauty in physics rather than mathematics, his viewpoint can be extended to encompass mathematics as well. In the subsequent paragraphs, we aim to demonstrate that the beauty encountered in mathematics is contextual, akin to the beauty of the racehorse in the earlier example. For instance, a mathematical proof is not beautiful in the same way a flower is; its beauty lies in its ability to serve a specific purpose — to offer insight into the realm of mathematics in a way that evokes fascination or admiration. 

Weinberg's example also serves to illustrate the view that our judgments of beauty in mathematics are grounded in facts rather than subjective preferences. Although the aesthetic experience and the associated feelings of fascination or admiration are inherently subjective, their evocation in mathematics is always tied to an objective property or relation between objects under consideration (compare this with how we subjectively experience warmth, this experience being caused by the objective fact of sitting next to a fireplace). The horse trainer’s judgment about the horse’s beauty is based on the properties of the horse that reflect its ability to win races, rather than on the trainer’s subjective preferences (such as a preference for a certain color or the length of the tail). In the same way, we argue that the property or set of properties evoking this particular type of aesthetic experience in mathematics is not relative to a mathematician's subjective preferences but is rather justifiable on objective grounds.  

But how exactly to define beauty and explain its connection to truth in mathematics? 
According to Weinberg's proposal, beauty is supposed to help a mathematician select ideas that can better explain the world of mathematics, and thus be closer to truth. 
Along these lines, we offer an account of beauty according to which:  

\textit{Beautiful mathematics is one in which a great deal of the entire mathematical world is reflected}.

Our definition admits of varying degrees of beauty depending on the intensity and extent of this reflection. In other words, the more an object reflects the totality of the mathematical world, the more beautiful it is. This criterion should apply not only to mathematical objects but to mathematical results as well. A beautiful result allows better and wider insight into the mathematical reality as a whole. It reveals connections between apparently non-related areas or problems, it is closely related to many other results and thus it helps us expand our view of mathematical reality. This explains a close connection between epistemic and aesthetic values in mathematics. If a result provides better and wider insight connecting non-related areas, it is likely to be more fruitful and explanatory potent. 

The focus of mathematics is not on understanding specific objects and facts, but rather on using objects and facts as lenses through which we can observe more general laws and the intricate structures of the mathematical universe as a whole. A piece of mathematics is considered beautiful inasmuch as it helps us attain that goal and drives our attention toward results and properties that offer a more general and deeper understanding of the subject. 
Some other views also connect mathematical beauty to some form of comprehension (see \cite{cellucci} and \cite{rota}). Rota centers his account of mathematical beauty on the concept of enlightenment, and describes beautiful mathematics as one that is enlightening, while Cellucci prefers the term understanding, as enlightenment might misleadingly suggest that appreciating mathematical beauty is instantaneous and doesn't presuppose prior knowledge of the subject.  

Our perspective ultimately draws inspiration from Plato and his idea about the role of beauty in acquiring knowledge of Forms. In the following chapter, we briefly discuss Plato's views on beauty and their connection to mathematics.  

\section{Plato on revealing the unity beneath diversity}

In Timaeus, Plato narrates how god created the universe using the most beautiful mathematical objects as a model. By a beautiful mathematical object, Plato means ``a straight line or a circle and resultant planes and solids [...] these things are not, as other things are, beautiful in a relative way, but are always beautiful in themselves'' (\cite{fileb}, p. 51). With shapes and numbers, god brought order into the chaos that used to rule. He gave the universe a perfectly spherical shape ``because there is no shape more perfect and none more similar to itself -- similarity being, in his opinion, incomparably superior to dissimilarity'' (\cite{timaj}, 33b). This refers to the fact that a perfect sphere is completely symmetrical around its center (it has infinitely many symmetrical transformations). Using other superiorly beautiful figures, god created fire, air, water, and earth. He made them ``as beautiful and as perfect as they could possibly be'' (\cite{timaj}, 53b). These figures are so-called Platonic solids: tetrahedron, octahedron, cube, and icosahedron. They are composed of two basic triangles. Both of them have a right angle, but one is isosceles, and the other is, according to Plato, the most beautiful of scalene triangles, that is, the half of an equilateral triangle. Besides these four beautiful solids, Plato mentions dodecahedron as the shape of the elements of ether.

According to Plato, the perfect sphere is aesthetically superior to the other geometric bodies because it is the most similar to itself, that is, the most symmetrical and proportionate. The beauty of the five solids can be explained similarly. The symmetry of these solids is the consequence of the fact that their faces are identical regular polygons, with the same number of them meeting in each of the vertices. These objects are nowadays known as \textit{regular polyhedra}. 

In Plato's philosophy, the experience of beauty, including that of mathematical objects, plays a significant role in the pursuit of ultimate knowledge which is for Plato the knowledge of Forms. Forms, according to Plato, are universals or concepts that manifest themselves in individual objects but exist independently of them. Plato believed that only Forms have real existence, being eternal and changeless, while all other concrete objects are merely their reflections. Hence, he considered the knowledge of Forms to be the ultimate and the only true knowledge. The soul, being immortal, has already come to know the Forms, and the process of learning is the process by which the soul recalls this knowledge.

In the ascent of the soul toward the perception of the Forms, an important role is ascribed to mathematics, the science that deals with measure, symmetry, and shape, in which beauty is contained.\footnote{As we use it here, the term `science' refers to any systematic and methodological study of reality.} Beautiful mathematical objects, such as those mentioned above, are especially likely to revive the recollection of the Forms. Namely, beautiful things remind the soul of beauty itself, which belongs to the world of Forms, and facilitate contemplation or recollection of this world. Through beauty ``we are brought to recognize and value the ideal properties of the Form of beauty and, indeed, of all the Forms'' (\cite{obdrzalek}, p. 205). Beauty has such an important role in prompting the recollection of the Forms because the Form of beauty can ``shine through its instances'' and show in them its ideal properties (\cite{obdrzalek}, p. 214). It helps the soul find \textit{the unity in diversity}, that is, recognize the Form of beauty in its different manifestations.

The motivation for our view on the beauty in mathematics can be found in different aspects of Plato's philosophy. Plato took beauty as something objective and closely related to the possibility of gaining real knowledge, which is a view we aim to develop here applying it to the beauty in mathematics. In this, we follow Plato's idea about the nature of real knowledge. Plato believed that this knowledge is not related to a particular object and its specific properties. Rather, it can be reached only when it is seen past particularities and when the aspects of objects that represent universal laws or ideas, i.e. Plato's Forms, are brought to light. These Forms unite objects that on the surface might seem very diverse.

In the next sections, we present our account of beauty in mathematics inspired by Plato's views on the role of beauty in revealing ``the unity beneath diversity''. We explain what this means in the context of mathematics and how particular features of beautiful mathematical objects, such as symmetry, follow from this account. Before that, we distinguish different types of mathematical beauty and present their examples.

\section{Three notions of beauty in mathematics}

We distinguish three types of mathematical beauty, each corresponding to a different category of entities that can be characterized as beautiful: 
\begin{itemize}
    \item[$\cdot$] the first one is the beauty of \textit{mathematical objects} and the structures they form;
    \item[$\cdot$] the second one is the beauty of \textit{theorems} that describe these objects; 
    \item[$\cdot$] the third type encompasses the beauty found in the \textit{proofs} of these theorems and the tools they employ, including axiomatizations, definitions, methods and techniques, etc.
\end{itemize}

The following metaphor describes how one can understand the relationship between the three categories of mathematical entities. The mathematical world, consisting of mathematical objects and the structures they form, resembles a landscape. Mathematical theorems are like photographs of the landscape, which can be captured from various angles and focused on distinct elements. On the other hand, proofs and other techniques used in reaching these theorems can be compared to camera lenses, which allow us to view the mathematical landscape from specific vantage points. The camera enables us to capture pictures of a landscape with varying information content and clarity. This depends on the specific part of the landscape we choose to focus on, our camera setup, and the lenses used.

According to our perspective, the beauty of any mathematical entity is connected to its place within the corresponding level of mathematics. Its beauty arises from its relationship with other entities at the same or different levels, rather than from any specific properties of the object itself. 

In the following chapters, we provide examples of the three types of beauty which we hope will further illuminate our perspective. We also try to determine the common characteristics that unify these notions into a general underlying concept of mathematical beauty. Although mathematical objects, theorems, and proofs have different types of aesthetic qualities and manifest different aspects of beauty, we believe that one of these notions is fundamental, as it serves as the underlying basis for the other two. 

\subsection{Beauty of mathematical objects}

The following examples show the difference between the aesthetic appeal of mathematical objects and other kinds of beauty we find in art and nature. 

\begin{figure}[htb]
    \centering
\begin{subfigure}{0.3\textwidth}
\centering
    \includegraphics[height=4.5cm]{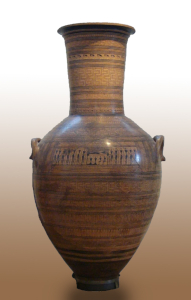}
    \caption{Vase}
    \label{fig:a}
\end{subfigure}%
\begin{subfigure}{0.3\textwidth}
\centering
    \includegraphics[height=4.5cm]{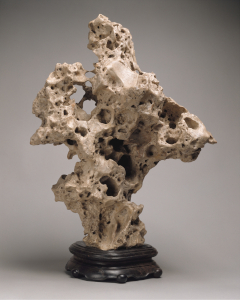}
    \caption{Gongshi}
    \label{fig:b}
\end{subfigure}
\begin{subfigure}{0.3\textwidth}
\centering
    \includegraphics[height=4.5cm]{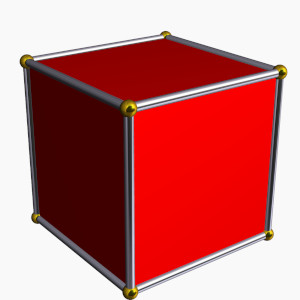}
    \caption{Cube}
    \label{fig:c}
\end{subfigure}%

    \caption{}
    \label{fig:objects}
\end{figure}

Let us consider three objects - (A) a painted Greek vase, (B) Gongshi\footnote{Gongshi, also known as scholar's rocks, are naturally occurring rocks which are traditionally appreciated by Chinese scholars.}, and (C) a cube - and ask ourselves about their aesthetic value. If an ordinary person, not particularly interested in mathematics, were asked which of them is (most) beautiful, the answer would probably be (A) or (B).

However, if a mathematician were asked which of the three has the greatest aesthetic value of the mathematical type, the answer would be (C).
To see more clearly why, let us examine what aesthetic values these three objects exemplify, and how they relate to mathematics.

(A) Painted Greek vase. This is clearly an art object, whose aesthetic value is assessed according to the criteria we use to assess the beauty of any artwork. However, the fact that the object is a vase, possessing that particular shape which is characterized by its distinct conic form, is an important aspect as well. If we were to remove the painting and abstract the shape by eliminating the handles, we would observe that the shape we are left with is very mathematical. 
As a mathematical object, it reflects the concepts of circle and exhibits perfect rotational symmetry. Therefore, the beauty of the vase as an artwork is intertwined with the beauty of the mathematical objects it takes its shape from.

(B) Gongshi. Here we have an object that occurs in nature -- a stone that is awkwardly asymmetrical and textured but also beautiful in its own right. Being an object of nature, the formation of this rock can be explained through laws of physics which are in turn expressed through mathematics. For example, the eroded surface of the rock can be modeled with certain dynamical systems; the structure of the rock can be explained by crystallography whose laws are governed by Euclidean geometry, etc. 
So much mathematics can be connected to the object, but none of that mathematics is specific to it nor is it inherent in it in a way the conic shape is inherent in (A). In other words, the beauty of this object has nothing to do with the mathematics used to describe it. Since no mathematical object or concept fundamentally underpins the object, we can say that it doesn't reflect any mathematics and hence it does not possess the beauty of mathematical type.

(C) Cube. There are very few things occurring in nature that resemble a cube and we might safely say that a cube is not an object of nature at all. And while cubic shapes can occur in artwork (for example, we could have had a painted cube instead of a painted vase), cube qua cube is no work of art. It is a mathematical object and it is considered a very beautiful one. But the aesthetic appeal of a cube is different from the aesthetic appeal of a vase or a gongshi. 
The conic shape of a vase, we saw, reflects the concept of rotational symmetry. The cube, similarly, reflects the notion of cubic symmetry. However, unlike the conic shape of the vase, the cube reflects so much more than its particular symmetry. Indeed, it sits at the confluence of what we would classify as many distinct branches of mathematics.  
In no particular order, one sees the cube as (1) convex polyhedron, (2) regular polytope, (3) zonohedron, (4) Wythoff, i.e. kaleidoscopic construct, (5) unit ball of a metric, (6) topological sphere, etc.
All these ``faces'' of the cube open a window to whole fields of mathematics. Moreover, by looking at all the ``faces'' together we glimpse at the unity of all those fields as well. This makes a cube an object of mathematical beauty par excellence. As already mentioned, our account allows for varying degrees of beauty in mathematics. So one can say that a cone, although a beautiful mathematical object, is however less beautiful than a cube since it reflects a smaller part of mathematical reality.

\medskip

In some of our examples, and across literature on mathematical beauty, a prominent place is given to the concept of symmetry. Our account explains why this is the case.
The symmetry is, naturally, often thought of as ``geometrical symmetry'', which has an important place in visual arts as well. However, symmetry is a much deeper concept, which we'll briefly clarify.
An invariant is a property of a mathematical object that remains unchanged when transformations of a certain type are applied to it.
For example, the degree of a polynomial is invariant under the linear change of variables; the area of a figure in the plane is invariant under isometries of the plane, etc. 
Especially important are the transformations under which an entire object is left unchanged, i.e. under which the object itself is invariant. These transformations are called symmetries.
The set of transformations that fix invariant forms a group. In particular, when the transformations are symmetries, we speak of the group of symmetries. 
Some mathematical objects share their group of symmetries. For example, the group of symmetries of the cube is $O_h$ ($S_4 \times S_2$). The same is the group of symmetries of the octahedron (indeed the standard name for the group in question is ``octahedral group''). Therefore, cube, octahedron, and group $O_h$ (a bona fide mathematical object in its own right) are directly and inseparably connected. In our somewhat poetic language, they are all reflected in any of them. The more objects that have a particular group of symmetries in common, the larger segment of mathematics each of them reflects, and the more mathematically beautiful they are.

According to our account, the concept of symmetry occupies such a significant role in discussions of mathematical beauty because it makes the connections between different mathematical objects particularly evident. Nevertheless, it is only one of many ways in which these connections occur. As we've seen in the case of a cube (and we'll see in other examples to follow), mathematical objects and areas are connected in other important ways. Each connection allows an object to reflect related objects and makes it more beautiful.  

The kind of beauty that consists in the ability to reflect different objects and areas is specific to mathematics and it isn't necessarily valued in art or nature. In these areas, what sets an object apart as particularly beautiful often has nothing to do with its connection to other things. Rather, what is valued is usually the object's particularity and uniqueness. The appeal of expressionistic paintings, for example, lies in distortions of reality which serve to express a unique and subjective vision of the artist. In this regard, mathematical beauty is very different from the beauty of artwork. It comprises a separate aesthetic category, which is bound and driven by a specific purpose of describing and understanding the abstract domain of mathematics.  

\subsection{Beauty of theorems}
In Hardy's words, mathematical theorems are ``notes of our observation'' of mathematical reality (\cite{hardy}, p. 124). Although their purpose is to provide us with insight into the world of mathematics, they can also be seen as something that possesses aesthetic value. But what makes a theorem aesthetically appealing and what is the relationship between the beauty of a mathematical theorem and the beauty of the object it depicts? 

In line with our metaphor (p. 7), a theorem allows us to paint a picture of the mathematical world. A picture can be beautiful because it depicts a beautiful part of reality, but also because it is made in a beautiful way. Similarly, theorems can be beautiful because they describe a beautiful object or a relation between beautiful objects. On the other hand, they can also be appreciated because of the light they throw on these objects and their relations. Beautiful theorems might allow us to observe them from a particularly suitable angle, or to recognize their particularly important features, thereby illuminating an aspect of mathematical reality that significantly improves our understanding of it. 

According to this account, the beauty that characterizes a theorem doesn't represent a new kind of beauty. Rather, it depends on the beauty of its content. A beautiful theorem allows us to recognize beautiful aspects of mathematical reality and it brings our attention to general and recurring features of mathematical objects that make them beautiful. The theorem, therefore, owes its beauty to the beauty of the object it describes. 

This view explains why mathematicians often say that there can be no beautiful theorems about trivial facts or uninteresting objects. 
It can also explain how the beauty of mathematical theorems, as we see it, can be related to scientifically valuable results. According to Poincare, those are the results that show similarities hidden under apparent discrepancies and ``unite elements long since known, but till then scattered and seemingly foreign to each other'' (\cite{poincare}, p. 30). They introduce order into apparent disorder and thus help us navigate physical or mathematical reality. 

\subsection{Beauty of proofs}
A mathematical proof is not a mere verification. Proofs are the means by which we come to comprehend mathematics and advance our knowledge. They too can be subjected to aesthetic assessment and are often judged as more or less clear, concise, simple, motivated, elegant, or beautiful. Unlike theorems, proofs appear to possess a twofold nature. While essentially dependent on the theorem's content (as suggested by Rota, who posited that there are likely no beautiful proofs of ugly theorems in \cite{rota}, p. 172), proofs also reflect our approach to understanding and manipulating this content. How do these two facets of proof interplay in endowing it with aesthetic properties? 

On one hand, the aesthetic appeal of a proof can stem from the method used. Techniques such as visualization or diagonalization, for instance, can make the proof remarkably simple and effective. However, such proofs are better qualified as \textit{elegant} than strictly speaking \textit{beautiful}. According to Rota, ``mathematical elegance has to do with the presentation of mathematics, and only tangentially does it relate to content'' (\cite{rota}, p. 178). In other words, elegance is an aesthetic property that has nothing to do with \textit{what} is proved but only with \textit{how}.\footnote{For an interesting analysis of elegance of proofs see \cite{harre}.}

On the other hand, some proofs are called beautiful because they reveal the essence of the theorem they prove. The beauty of such proofs hinges on the theorem's content as well as the proof's capacity to elucidate it and to reveal \textit{why} the theorem holds (cf. \cite{mordell}). The proofs that can help us explain the \textit{why} of the theorems, so to speak, are also called \textit{deep}. They are often general, in the sense that they contain ideas that can be applied to a family of similar problems. Such ideas, ``constituent in many mathematical constructs'' and ``used in the proof of theorems of many different kinds'', Hardy calls \textit{significant} (\cite{hardy}, pp. 103-109). Deep proofs are based on some fundamental properties of the mathematical world, not confined to a particular object or objects of a particular kind. 

The beauty of mathematical proof thus rests on two pillars -- first, its effectiveness in verifying the proposition in question (its elegance); and second, its connection to general mathematical concepts (its depth). We will try to illuminate the distinction by presenting some examples of proofs. We engage in their comparative analysis, aiming to determine which one exhibits a greater degree of beauty and why. Throughout this examination, we delineate the mathematical facets contributing to their aesthetic appeal. This involves assessing whether the beauty we find in these proofs falls into our category of elegance or depth. 

\begin{proof}[(A) Proof of Pythagoras Theorem]
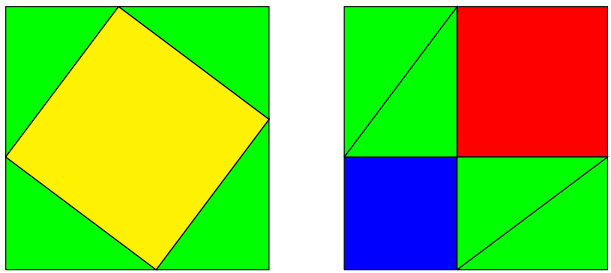
\begin{figure}[H]
  \centering
  \input{pythagoras}
  \caption{(A) Proof of Pythagoras Theorem.}
  \label{fig:pythagoras}
\end{figure}
To see that the two pictures contain proof of Pythagoras theorem, one only needs to realize that the yellow square in the first picture can be decomposed into two squares from the second one. This is accomplished by rearranging the four triangles. Since the yellow square represents the square over the hypothenuse $c$, the red square is a square over the longer of the remaining sides $a$, and the blue square over the shorter one $b$, this shows that $c^2=a^2+b^2$.
\end{proof}

\begin{proof}[(A') Another proof of Pythagoras Theorem]
\begin{figure}[ht!]
\centering
\includegraphics[width=50mm]{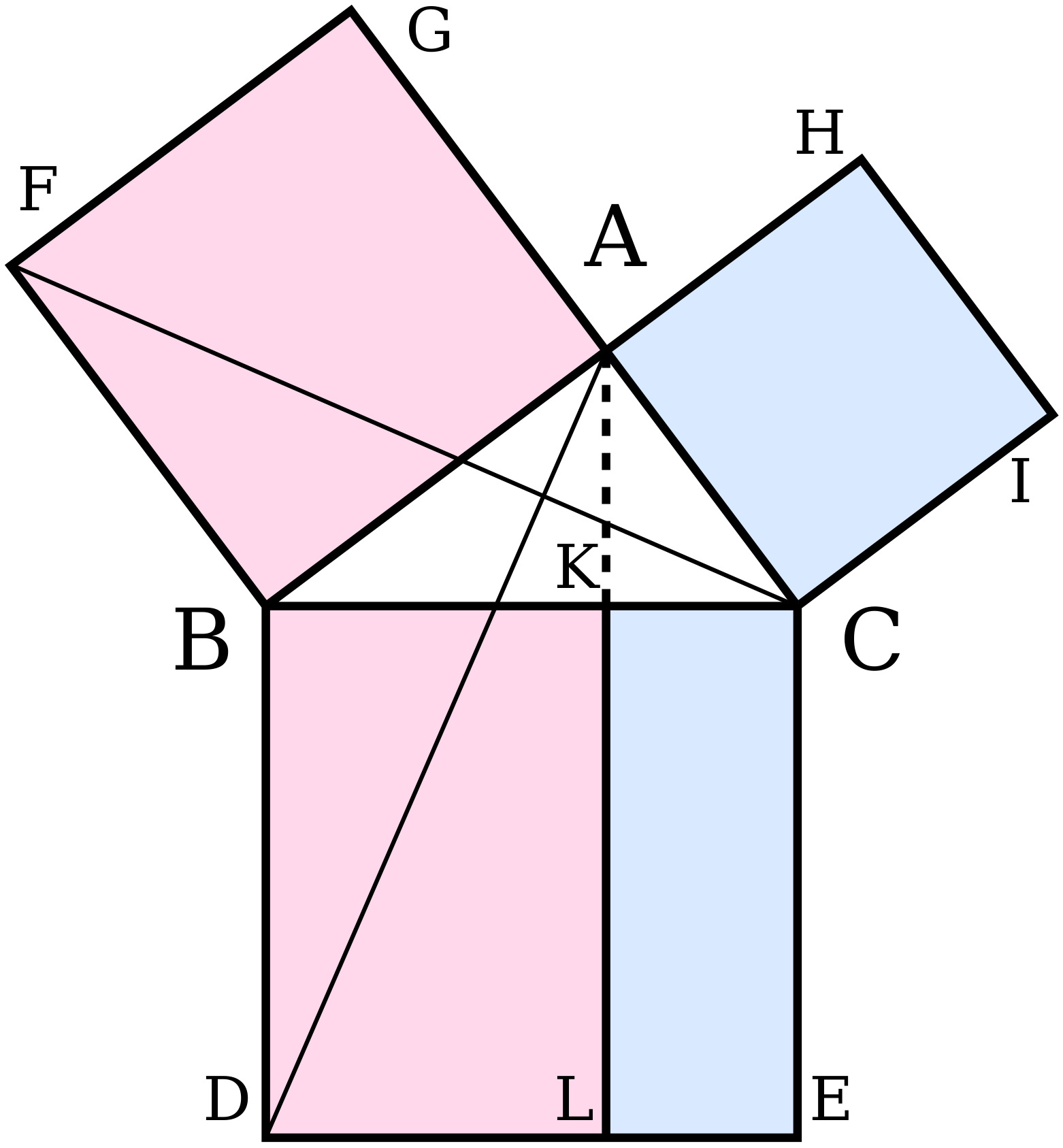}
\end{figure}
Here, the idea is to show that the areas of the same color are equal, from which Pythagoras theorem easily follows. That pink-colored areas are equal is shown by constructing the two triangles, $ABD$ and $FBC$. It is first shown that the triangle $ABD$ is half of the rectangle $BDLK$ using the lemma according to which the area of a triangle is half the area of any parallelogram on the same base ($BD$) and having the same altitude ($AK$). For the same reason, the triangle $BFC$ is half of the square $ABFG$. Next, it is shown that the two triangles are congruent because they share two sides and the angle between them. The equality of the blue areas is shown similarly.   
\end{proof}

Most mathematicians would agree that (A) is a beautiful proof. It is witty and a bit unexpected, mostly because it can be expressed entirely in two neat pictures. What is also pleasing is the ease with which we come to realize that the theorem holds just by looking at the pictures. However, the aesthetic quality of the proof derives in great part precisely from those pictures that are in fact perceived as art (they resemble an abstract art piece). Setting aside its artistic appeal, one is left with an effective and elegant argument that, in the course of verifying the theorem in question, presents the important ideas of congruence and equidecomposability. These ideas make the proof in question deep. 

The proof (A') of the same theorem is surely less beautiful. It uses an auxiliary construction that helps to show that the theorem holds but it does not reveal \textit{why}. It is made purposely to arrive at the proof in question and lacks further motivation. The truths it relies on concern the equality of two figures and the means of establishing it. They cannot be used for proving similar theorems (such as the one saying that the square over the side opposite the obtuse angle of the triangle is strictly greater than the sum of the squares over the two other sides). The proof doesn't contain deep ideas that reappear in different places and that can be used in proving other, even similar, theorems. That is why it is not particularly deep. It is also much less elegant than the proof (A), i.e. it is longer and more difficult to follow (for a similar analysis of the aesthetic worth of this proof and its comparison to another proof of Pythagoras theorem, see \cite{cellucci}, pp. 347-348). 

\begin{proof}[(B) Proof of the infinity of primes]
Let $\mathbb{P} \subset \mathbb{N}$ be the set of all primes. Suppose that $\mathbb{P}$ is finite, i.e., $\mathbb{P}=\{p_1,\ldots,p_n\}$. Let $N=p_1 \cdots p_n +1$. Clearly, $N>1$ and is indivisible by any $p \in \mathbb{P}$. Therefore, there must exist a prime $P$ (possibly $N$ itself) that is not in $\mathbb{P}$, which contradicts the assumption that $\mathbb{P}$ is the set of all primes. Hence, $\mathbb{P}$ is infinite.
\end{proof}

In comparison to the proof in (A), the proof in (B) has no artistic elements that need to be set aside -- its content is purely mathematical. It too gives an effective and elegant argument (perhaps slightly less elegant than (A)s due to its slightly lengthier form) and presents an important idea of ``proof by contradiction'' and an archetype of a constructive method of refuting the assumed finiteness of a set. Given how widespread this idea is (much more so than the strictly geometrical ideas of congruence and equidecomposability shown in (A)) we can say with confidence that (B) is more beautiful than (A).

\begin{proof}[(C) Proof of the logical truth $((p\rightarrow q)\land p)\rightarrow q$]
Consider first its proof in a Hilbert-style, axiomatic formal system:
\begin{scprooftree}{0.6}
\AxiomC{$(((p\rightarrow q)\land p)\rightarrow (p\rightarrow q))\rightarrow ((((p \rightarrow q)\land p)\rightarrow p)\rightarrow (((p\rightarrow q)\land p)\rightarrow q))$}
\AxiomC{$((p\rightarrow q)\land p)\rightarrow (p\rightarrow q)$}
\BinaryInfC{$(((p \rightarrow q)\land p)\rightarrow p)\rightarrow (((p\rightarrow q)\land p)\rightarrow q)$}
\AxiomC{$((p \rightarrow q)\land p)\rightarrow p$}
\BinaryInfC{$((p \rightarrow q)\land p)\rightarrow q$}
\end{scprooftree}
\end{proof}

\begin{proof} [(C') Another proof of $((p\rightarrow q)\land p)\rightarrow q$]
This proof is made inside a natural deduction formal system:
\begin{prooftree}
\AxiomC{$[(p\rightarrow q)\land p]$}
\UnaryInfC{$p\rightarrow q$}
\AxiomC{$[(p\rightarrow q)\land p$]}
\UnaryInfC{$p$}
\BinaryInfC{$q$}
\UnaryInfC{$((p\rightarrow q)\land p)\rightarrow q$}
\end{prooftree}
\end{proof}

Even though here we deal with formal proofs of logical truth, we can still discern properties similar to those found in mathematical proofs, that make some of them more beautiful than others. 

The proof in (C) is a result of a clever use of available axioms. All we need to do is find a way to instantiate those axioms that allow us to infer the formula to be proven by available rules of inference. It shows how $((p\rightarrow q)\land p)\rightarrow q$ follows from the distributivity of implication and the axiom standing for the elimination of conjunction, by repeated application of modus ponens. Since these axioms are logical truths, and modus ponens is a sound rule of inference, we can be sure that $((p\rightarrow q)\land p)\rightarrow q$ is a logical truth as well. However, the idea of the proof remains unclear. How did we know which axioms to use in the proof, and what is the role of distributivity in making $((p\rightarrow q)\land p)\rightarrow q$ true? Are the logical truths we used in the proof more fundamental than the truth we proved? The proof lacks insight into the meaning of the proven formula and fails to introduce any significant ideas that could be used in dealing with other formulas of the same kind. 

Axiomatic logical systems abound with proofs that seem to be unmotivated by the formula to be proven and lack a clear, potentially generalizable, idea. This raises doubts about the suitability of these systems for rigorous study of deduction. 

Alternative formal frameworks, such as natural-deduction systems, offer a more intuitive and coherent approach to exploring logical truths and deductions. The proof in (C') uses general ideas concerning deduction that are closely related to the meaning of implication. It is based on a simple and clear idea of how an implication is to be proven: by assuming its antecedent and deducing the consequent from it. If we manage to do this, we have proven the implication and our proof no longer depends on the assumption that its antecedent holds. This assumption can thus be canceled. The proof in (C') thus unveils a crucial aspect of implication: an implicative statement can be interpreted as asserting the potential deduction of its consequent from its antecedent. Not only are the ideas used in the proof simple and natural, but they are also very general in the sense that they can be used whenever we want to prove an implication. They encode a very natural way of reasoning about implicative statements that is used all over logic and mathematics. This makes the proof in (C') deeper than the one in (C). Moreover, its brevity, directness, and ease of comprehension make it elegant as well.

The ideas contained in the natural-deduction proof are, in a sense, explicated by its following version in sequent calculus:

\begin{prooftree}
\def\fCenter{\ \vdash\ }
    \AxiomC{$p \fCenter p$}
    \AxiomC{$q \fCenter q$}
    \BinaryInfC{$p\rightarrow q, p \fCenter q$}
    \UnaryInfC{$(p\rightarrow q) \land p \fCenter q$}
    \UnaryInfC{$\fCenter ((p\rightarrow q) \land p)\rightarrow q$}
\end{prooftree}

The turnstile can be read as `the formula on the right-hand side is deducible from the formulas on the left-hand side'. The inference steps then show what can be inferred about the deducibility relation between some other formulas. Read in this way, the second row describes the rule of modus ponens. The inference to the third row corresponds to the conjunction elimination rule as it is used in natural-deduction proofs: it shows that if a formula is deducible from the two formulas, then it is also deducible from their conjunction. The last step describes the main idea of the proof of implication which is that we prove an implication if we prove that its consequent is deducible from its antecedent. 

The sequent calculus makes the concept of deducibility central. It allows us to explore the connections between deducibility of different formulas from the same or different sets of formulas. It thus represents an important part of the so-called \textit{general proof theory} that deals with proofs in their own right, and attempts to understand their nature and mutual relations. This is one of the central goals of logic that could hardly be approached using axiomatic formal systems. So, the proofs from the natural-deduction or sequent calculus are not only more elegant, but they also allow us to approach the subject of logic in a more direct and revealing way (cf. \cite{DA}, pp. 32-36). They reflect other proofs using similar ideas. Moreover, considered in this way, deductions can be seen to build structures that reflect those found in algebra, topology, set theory, etc. This is what makes these proofs and their study beautiful according to our account.

\section{The unity of mathematical beauty}
At the beginning of the paper, we distinguished three kinds of beauty found in mathematics, each pertaining to different kinds of mathematical entities -- mathematical objects, theorems, or proofs. However, the subsequent discussion revealed the characteristics the three notions of beauty have in common. 

According to our account, each notion of beauty in mathematics has to do with the possibility of obtaining a wider and more informative insight into the mathematical world. Mathematical objects are beautiful if they reflect general features of the mathematical world that are found across different areas of mathematics. They are the more beautiful the greater share of the mathematical world they reflect. Theorems have a role to improve and widen our understanding of the mathematical world and their aesthetic appeal depends on how well they accomplish this task. It is best accomplished by the theorems that direct our attention to the very general aspects of mathematical objects shared with many others. But these are exactly those aspects that make the objects beautiful. So, the beauty of theorems comes from the beauty of the objects they describe. Finally, a mathematical proof is the more beautiful (i.e. deep) the more and the deeper concepts it uses. The possibility of using the same concepts and ideas in proofs regarding mathematical objects from different areas depends on connections between the areas that are manifested in the beauty of these objects.

The suggested perspective entails that the beauty of mathematical objects is the most fundamental type of beauty in mathematics. Put metaphorically, the beauty in mathematics is most akin to the beauty of a landscape. Its pictures, in the form of theorems, are beautiful to the extent to which they manage to represent faithfully its beauty, and the methods of making the pictures, i.e. of proving theorems, are beautiful to the extent to which they make this possible. 

Since the beauty of mathematical objects depends on properties that do not characterize beautiful objects from different areas, mathematical beauty, in general, is to be considered a \textit{sui generis} aesthetic category.

\begin{quote}
\textit{Acknowledgements: We would like to thank Milo\v{s} Ad\v{z}i\'{c} for his help on improving the text with his comments and suggestions.}
\end{quote}

\end{document}

%% file: pythagoras.tex
\begin{tikzpicture}[scale=0.5]
  \draw[color=black, fill=green] (-8,0) -- (-4,0) -- (-8,3) -- (-8,0);
  \draw[color=black, fill=green] (-4,0) -- (-1,0) -- (-1,4) -- (-4,0);
  \draw[color=black, fill=green] (-1,4) -- (-1,7) -- (-5,7) -- (-1,4);
  \draw[color=black, fill=green] (-5,7) -- (-8,7) -- (-8,3) -- (-5,7);
  \draw[color=black, fill=yellow] (-4,0) -- (-1,4) -- (-5,7) -- (-8,3) -- (-4,0);

  \draw[color=black, fill=blue] (1,0) -- (4,0) -- (4,3) -- (1,3) -- (1,0);
  \draw[color=black, fill=red] (4,3) -- (8,3) -- (8,7) -- (4,7) -- (4,3);
  \draw[color=black, fill=green] (4,0) -- (8,0) -- (8,3) -- (4,0);
  \draw[color=black, fill=green] (4,0) -- (8,3) -- (4,3) -- (4,0);
  \draw[color=black, fill=green] (1,3) -- (4,3) -- (4,7) -- (1,3);
  \draw[color=black, fill=green] (1,3) -- (4,7) -- (1,7) -- (1,3);
\end{tikzpicture}

%% file: main.bbl
\begin{thebibliography}{9}

\bibitem{rota}
Rota, G.C. (1997) ``The Phenomenology of Mathematical Beauty'', \textit{Synthese} vol. 111, pp. 171-182.

\bibitem{mcallister1}
McAllister, J. (1996) \textit{Beauty and revolution in science}, Ithaca: Cornell University Press.

\bibitem{thomas}
Thomas, R.S.D. (2016) ``Beauty is not all there is to Aesthetics in Mathematics'', \textit{Philosophia Mathematica} vol. 25, pp. 116-127.  

\bibitem{gerwen}
van Gerwen, R. (2011) ``Mathematical beauty and perceptual presence'', \textit{Philosophical Investigations} vol. 34, pp. 249-267.

\bibitem{cellucci}
Cellucci, C. (2015) ``Mathematical Beauty, Understanding, and Discovery'', \textit{Foundations of Science} vol. 20, pp. 339-355. 

\bibitem{todd}
Todd, C. (2008) ``Unmasking the Truth Beneath the Beauty: Why the
Supposed Aesthetic Judgements Made in Science May Not Be Aesthetic at All'', \textit{International Studies in the Philosophy of Science} vol. 22, pp. 61-79. 

\bibitem{harre}
Harr\'{e}, R. (1958) ``Quasi-aesthetic appraisals'', \textit{Philosophy} vol. 125, pp. 132-137. 

\bibitem{DA}
Do\v{s}en, K., and Ad\v{z}i\'{c}, M. (2017) ``G{\"o}del on Deduction'', \textit{Studia Logica} vol. 107, pp. 31-51.

\bibitem{vajnberg}
Weinberg, S. (1994) \textit{Dreams of a final theory}, New York: Vintage Books.

\bibitem{mur}
Moore, G. E. (1922) \textit{Principia Ethica}, Cambridge: Cambridge University Press.

\bibitem{fileb}
Plato (1975) \textit{Philebus}, Oxford: Clarendon Press.

\bibitem{timaj}
Plato (2008) \textit{Timaeus and Critias}, Oxford: Oxford University Press.

\bibitem{obdrzalek}
Obdrzalek, S. (2022) ``Why Er\={o}s'', in \textit{The Cambridge Companion to Plato}, second edition, Cambridge: Cambridge University Press, pp. 202-232.

\bibitem{hardy}
Hardy, G.H. (2012) \textit{A Mathematician's Apology}, Cambridge: Cambridge University Press.

\bibitem{poincare}
Poincar\'e, H. (1914) \textit{Science and Method}, London: Thomas Nelson and Sons.

\bibitem{mordell}
Mordell, L.J. (1959) \textit{Reflections of a Mathematician}. Canadian Mathematical Congress. 


\end{thebibliography}
